
\font\magbf=cmbx10 scaled\magstep2


\def\mbox#1{\leavevmode\hbox{#1}}
\def\frac#1#2{{#1\over #2}} 
 

\def\zz{\mbox{${\bf Z}$}}	
\def\qq{\mbox{${\bf Q}$}}	
\def\cc{\mbox{${\bf C}$}} 

\def\calc{\mbox{${\cal C}$}}
\def\calg{\mbox{${\cal G}$}} 

\count35=0
\count36=0
\count38=0


\def\section#1{\advance\count35 by 1
\vskip4ex\noindent{\magbf\number\count35 \ \ \ #1}
\count36=0\count38=0 \hfill\vskip2ex}

\def\subsection{\vskip1ex\advance\count36 by 1
{\noindent{\bf\number\count35.\number\count36}\ }} 

\def\today{\ifcase\month\or January\or February\or March\or April\or May\or
June\or July\or August\or September\or October\or November\or December\fi
\space\number\day, \number\year}

\def\beginth#1{\vskip3ex\advance\count36 by 1 {\noindent{\bf
\number\count35 .\number\count36} \ \ {\bf #1.} }\begingroup\sl} 

\def\endth{\endgroup\hfill\vskip2ex}

\def\defn{\vskip3ex\advance\count36 by 1 {\noindent{\bf \number\count35
.\number\count36} \ \ {\bf Definition.} }}

\centerline{\magbf   Obstructing Four-Torsion in the Classical}
\vskip2ex
\centerline{\magbf   Knot Concordance Group} 

\vskip3ex
\centerline{Charles Livingston and Swatee Naik}

\centerline{February 15, 1998}

\vskip7ex Suppose that $K$ is a knot in $S^3$ with 2-fold branched cover
$M_K$.  Our main result is the following.
\beginth{Theorem}  If $|H_1(M_K)| = pm$ with $p$ a prime congruent to 3 
mod 4 and gcd(p, m) = 1, then $K$ is of infinite order in the classical
knot concordance group, 
$\calc_1$.\endth

Our interest in this result is its application to the study of 4-torsion 
in the concordance group. There are 11 prime knots of 10 or fewer
crossings, beginning with the knot $7_7$, that represent elements of order
4 in the algebraic concordance group.  A simple calculation using this
theorem yields:

\beginth{Corollary} No prime knot with fewer than 11 crossings represents
an element of order 4 in $\calc_1$.\endth

Of greater interest than obstructing individual knots from being of order
4 is that the obstruction depends only on an abelian invariant of the
knot.  Hence corollaries like the next one concerning the Alexander
polynomial of a knot, $\Delta_K(t)$, follow readily.

\beginth{Corollary} If $\Delta_K(t) = 5t^2 - 11t + 5$  then $K$ is of
infinite order in $\calc_1$.\endth

By way of contrast, according to Levine [L2], {\bf every} knot with
$\Delta_K(t) = 5t^2 - 11t + 5$  is of order 4 in  $\calc_{2k-1}$ for $k >
1$.

\section{Introduction}

In his classification of the knot concordance groups, Levine [L1] defined
the algebraic concordance groups,
$\calg_{\pm}$, of Witt classes of Seifert matrices and a homomorphism from
the odd-dimensional knot concordance groups 
$\calc_{4n \pm 1}$ to $\calg_{\pm}$.  The homomorphism is induced by the
function that assigns to a knot an associated Seifert matrix: it is an
isomorphism on
$\calc_k, k \ge 5$; on $\calc_3$ it is injective, onto an index 2 subgroup
in smooth category and surjective in the topological locally flat setting;
for $k = 1$ it is surjective.  However,  Casson and Gordon [CG1,CG2]
proved that on $\calc_1$ the kernel is nontrivial.  (Casson and Gordon's
original work applied  in the smooth setting, but their results are
 now known to hold in the topological locally flat setting as well, a fact
that follows from the existence of normal bundles in topological
4-manifolds, [FQ]. Similarly, our work applies in both categories.) 
Later, Jiang [J] extended Casson and Gordon's work to prove that the
kernel of Levine's homomorphism is infinitely generated.

Levine [L2] also proved that $\calg_{\pm}$ is isomorphic to an infinite
direct sum,
$\calg_{\pm} \cong \zz^\infty \oplus\zz_2^\infty \oplus\zz_4^\infty$. 
There is 2-torsion in 
$\calc_1$  arising from amphicheiral knots, but beyond this little is known
concerning torsion in
$\calc_1$.  Fox and Milnor [FM], in the paper in which knot concordance is
defined, made this observation concerning amphicheiral knots and asked if
there is torsion of any order other than 2.  This question reappears as
problem 1.32 of [K1, K2].  In a different direction, in  1977 Gordon [G]
(see also [K2], problem 1.94) asked whether every order 2 class in
$\calc_1$ is represented by an amphicheiral knot; as of yet the only result
bearing on this question is the observation that in higher dimensions the
answer is no, and in dimension 3  there are order 2 classes in
$\calg_{-}$ that cannot be represented by order 2 knots in $\calc_1$ [CM].
   
In this paper we will use Casson-Gordon invariants to derive results
concerning 4-torsion in the classical knot concordance group.  That
Casson-Gordon invariants can be used to show that an individual knot that
is of algebraic order 4  is  not of order 4 in concordance is not
surprising, though the examples presented here are the first; that the
method applies to the knot $7_7$ is pleasing in that this is the first
knot identified by Levine [L2] as a candidate to be of  order 4. It is
surprising  that the obstructions we find depend only on the abelian
invariants of the knot, such as the homology of the cover and the
Alexander polynomial. In contrast to this, if a Seifert form represents 0
in
$\calg_-$ then there is a slice knot having that Seifert form; as a
consequence, if a knot is algebraically slice then there is a slice knot
having identical abelian invariants. A related result states that any
integral polynomial
$\Delta(t)$ with
$\Delta(1) =
\pm1$ is an Alexander polynomial of some slice knot if and only if
$\Delta(t) = f(t)f(t^{-1})$ for some polynomial $f$ [T]. 

\vskip3ex
\noindent{\bf Outline} In the next section we review the definition of
concordance and of Casson-Gordon invariants.  Section 3 presents the
properties of the 3-dimensional bordism groups $\Omega_3(\zz_p)$ and 
$\Omega_3(\zz_p
\oplus
\zz)$ and Section 4 describes how bordism invariants can be extracted from
Casson-Gordon invariants.  An implication of this work is Theorem 4.3,
stating that if
$M^3$ is a 3-manifold with 
$|H_1(M_K)| = pm$ where $p$ is an odd prime and gcd($p, m$) = 1, then any 
Casson-Gordon invariant corresponding to a surjective $\zz_p$ character on
$M^3$ is nontrivial.  Section 5 reviews the use of Casson-Gordon invariants
to obstruct slicing.

The proof of Theorem 0.1 is  fairly technical, so we begin with two special
cases.  In Section 6 we prove the theorem for the prime 3.  This case is
especially easy since there is essentially only one $\zz_3$ character on
$M_K$ that must be considered.  Section 7 presents a restricted form of
Theorem 0.1, considering the case $p = 7$ but only proving that $K$ is not
of order 4.  The point of isolating this case is to indicate how one can
deal with the presence  of more than one essential character.  Finally, in
Section 8 the full proof of Theorem 0.1 is presented.  Section 9 discusses
the corollaries.

\section{Concordance and Casson-Gordon Invariants}

All our work holds in both the smooth and the topological locally flat
categories.  Throughout this paper $p$ will denote a fixed odd prime. 
Homology groups are always taken with $\zz$ coefficients unless
specifically noted otherwise.  Manifolds are all orientable and compact. 

\subsection{\bf The knot concordance group} 

 A knot
$K$ in $S^3$ is called {\sl slice} if $(S^3,K) = \partial(B^4,D)$ for some
embedded 2-disk in $B^4$.  Knots $K_1$ and $K_2$ are called {\sl
concordant} if $K_1\ \#\ \mbox{--}K_2$ is slice, where $\mbox{--}K$
represents the mirror image of
$K$. The set of concordance classes of knots forms an abelian group under
connected sum, denoted $\calc_1$.

\subsection{\bf Casson-Gordon invariants}  

Let $(M^3,\bar\chi)$ be a closed 3-manifold with a homomorphism
$\bar\chi\colon  H_1(M^3) \rightarrow
\zz_p
\oplus \zz$.    The bordism group $\Omega_3(\zz_p \oplus
\zz)
\cong
\zz_p$ (see Section 3 below), so $p(M^3,\bar\chi) =
\partial(W^4,\bar\phi)$, where
$W^4$ is a 4-manifold and $\bar\phi\colon  H_1(W^4) \rightarrow \zz_p
\oplus \zz$.  

Let $t(W^3,\bar\phi) \in L_0(\cc(t))$ denote the intersection pairing on
$H_2(W^4,\cc(t))$, where the field coefficients are twisted by the $\zz$
action given by multiplication by $t$ and the $\zz_p$ action  given by
multiplication by
$e^{2\pi i / p}$.  This pairing is viewed as an element in the Witt group
of nonsingular hermitian forms on finite dimensional $\cc(t)$ vector
spaces.  (In the case that the intersection form is singular, one must
first mod out by the radical of the form to achieve a nonsingular pairing.)

The invariant $\tau$ is defined by $\tau(M^3,\bar\chi) = {1 \over
p}(t(W^4,\bar\phi) - t_0(W^4))\in L_0(\cc(t)) \otimes \qq $, where $t_0$
is the class represented by the standard intersection form on
$H_2(W^4,\cc)$.

If $K$ is a knot in $S^3$, $M_K$ its 2-fold branched cover, and $\chi$ a
character from $H_1(M_K)$ to $\zz_p$, then there is a naturally induced
character 
$\bar\chi\colon  H_1(M_{K,0}) \rightarrow \zz_p \oplus \zz$, where
$M_{K,0}$ is the 3-manifold obtained from $M_K$ by performing 0-surgery on
the lift of $K$.  This follows from the fact that $H_1(M_{K,0})$ naturally
splits as  $H_1(M_K)
\oplus \zz$ with the generator of the $\zz$ factor given by the meridian of
the lift of $K$.  Hence, $\bar\chi$ is defined by mapping the meridian to
$(0,1) \in \zz_p \oplus \zz$.  
\beginth{Definition} The {\sl Casson-Gordon} invariant $\tau$ is defined by
$\tau(K,\chi) =
\tau(M_{K,0},\bar\chi)
\in  L_0(\cc(t)) \otimes \qq$.\endth  

An associated  signature invariant is defined as follows.  For a class in
$L_0(\cc(t))$ the signature is defined by evaluating a representative of
the class at a unit complex number and taking the limit of the signature
of the resulting complex valued form as the unit complex number approaches
$1$.  This map induces a homomorphism $\sigma\colon L_0(\cc(t)) \otimes
\qq \rightarrow \qq$. 

\beginth{Definition} The Casson-Gordon signature  invariant $\sigma$ is
defined by
$\sigma(M^3,\bar\chi) = 
\sigma(\tau(M^3,\bar\chi))$, and  $\sigma(K,\chi) = 
\sigma(\tau(K,\chi))$. \endth 

\noindent (Note, in [CG1] this is denoted
$\sigma_1\tau(K,\chi)$.) 

\vskip2ex
\subsection{\bf Additivity}

Given a knot $K = K_1\ \#\ K_2$, we have $M_K = M_{K_1}\ \#\ M_{K_2} $ and
any
$\zz_p$-valued character $\chi$ on $H_1(M_K)$ can be written as $\chi_1
\oplus
\chi_2$.  A key result of Gilmer [Gi] is that in this situation
$\tau(K,\chi) =
\tau(K_1,\chi_1) + \tau(K_2,\chi_2)$.

\section{Bordism Results: the groups $\Omega_3(\zz_p)$   and
$\Omega_3(\zz_p \oplus \zz)$}

The properties of the 3-dimensional bordism groups follow most easily from
a
 consequence of the bordism spectral sequence: $\Omega_3(G) = H_3(G)$.  Of
course they follow as well from more general bordism theory; [CF] is a good
reference.

First we have that $\Omega_3(\zz_p) \cong \zz_p$.  The map to $\zz_p$ is
given as follows.  For a pair $(M^3,\chi)$ with $\chi\colon H_1(M^3)
\rightarrow
\zz_p$ we view $\chi \in H^1(M^3,\zz_p)$.  The quantity $\chi \cdot
b(\chi)([M^3]_p)$ is the desired element in $\zz_p$.  Here $b$ represents
the Bockstein,
$b\colon H^1(M^3,\zz_p) \rightarrow H^2(M^3,\zz_p)$, the product is the cup
product, and $[M]_p$ is the $\zz_p$ reduction of the fundamental class of
$M$.

A useful alternative definition of the isomorphism is given using the
linking form, $\beta\colon \mbox{torsion}(H_1(M^3)) \times
\mbox{torsion}(H_1(M^3))
\rightarrow \qq/\zz$.  The restriction of $\chi$ to torsion($H_1(M^3)$) is
given by linking with some element $x \in $ torsion($H_1(M^3)$); that is,
$\chi(y) = \beta(x,y)$ for all $y \in \mbox{torsion}(H_1(M^3))$.  The 
self-linking of $x, \beta(x,x)$, is in $\qq/\zz$, but since it is
$p$-torsion, it can be viewed as an element in $\zz_p$.

We have the following result:

\beginth{Theorem} If $|H_1(M^3)| = pm$ with $p$ and $m$ relatively prime,
then for any nontrivial $\zz_p$ character, $[M^3,\chi]$ is nonzero in
$\Omega_3(\zz_p)$.\endth

\noindent{\bf Proof.}  For such a manifold the Bockstein is an
isomorphism.  By Poincar\'e duality the cup product is nontrivial in
$H^3(M^3,\zz_p)$.  Hence, $\chi \cdot b(\chi)$ is nontrivial in
$H^3(M^3,\zz_p)$ and the result follows.

\vskip2ex

We also will be using the fact that  map $\Omega_3(\zz_p)
\rightarrow \Omega_3(\zz_p \oplus \zz)$ induced by inclusion is an
isomorphism, with inverse induced by projection.  This follows from either
the Kunneth formula on homology or Kunneth results on bordism.  Again, see
[CF].

\section{Casson-Gordon Invariants as Bordism Invariants}

If $(M^3_1,\bar\chi_1)$ and $(M^3_2,\bar\chi_2)$ represent the same
element in the bordism group $\Omega_3(\zz_p \oplus \zz)$, then
$\tau(M_1^3,\bar\chi_1)$ and
$\tau(M_2^3,\bar\chi_2)$ differ by an element in $L_0(\cc(t))$.  Hence the
difference $\sigma(M^3_1,\bar\chi_1) - \sigma(M^3_2,\bar\chi_2)$ is an
integer. It follows that $\sigma$ defines a homomorphism $\sigma'\colon 
\Omega_3(\zz_p
\oplus
\zz)
\rightarrow \qq/\zz$.  This homomorphism takes values in
$(({1 \over p})\zz )/ \zz$, and so can be viewed as a
homomorphism
$\sigma_p\colon \Omega_3(\zz_p \oplus \zz) \rightarrow \zz_p$.  The induced
homomorphism on $\Omega_3(\zz_p)$ will also be denoted $\sigma_p$.

\beginth{Theorem} For $p$ odd, $\sigma_p$ is an isomorphism.\endth
\noindent{\bf Proof.}  A calculation of [CG1] shows that for the lens space
$L(p,1)$ given by $p$-surgery on the unknot in $S^3$ with $\zz_p$ character
taking value 1 on the meridian, $\sigma_p = 2$.  Since $p$ is odd, the
result follows.

\vskip2ex

\beginth{Theorem}  For a knot $K$ and $\zz_p$ character $\chi$ on
$H_1(M_K)$,
$\sigma_p(M_{K,0},\bar\chi) = \sigma_p(M_{K},\chi)$. Equivalently, 
$p\sigma(K,\chi) = \sigma_p(M_K,\chi)$ mod  $p$. \endth

\noindent{\bf Proof.} Since $\Omega_3(\zz_p \oplus
\zz) \rightarrow \Omega_3(\zz_p)$ induces an isomorphism, it follows that
$\sigma_p(M_{K,0},\bar\chi) =
\sigma_p(M_{K,0},\chi \oplus (0))$, where $(0)$ represents the trivial
character to $\zz$.  A $\zz_p$-bordism from $(M_{K,0},\chi \oplus (0))$ to
$(M_{K},\chi )$ is constructed from $M_K \times [0,1]$ by adding a 2-handle
with 0 framing to $M_K \times \{1\}$.  Note that $\chi$ extends over this
bordism since the 2-handle is added along a null homologous curve (the
lift of
$K$), and provides the desired $\zz_p$-bordism.\vskip2ex

We can now make one of our key observations.

\beginth{Theorem}  If $H_1(M_K) = pm$ with gcd$(p,m) = 1$ and $p$ odd,
then for any nontrivial $\chi\colon  H_1(M_K) \rightarrow \zz_p$,
$\sigma(K,\chi) 
\ne 0.$ \endth

\noindent{\bf Proof.} By Theorem 4.2, the nontriviality of $\sigma(K,\chi)$
will follow from that of $\sigma_p(M_K,\chi)$.  Theorem 3.1 implies that
$[M^3,\chi]$ is nontrivial in  $\Omega_3(\zz_p)$.  So Theorem 4.1 gives the
desired result.

\vskip2ex

\section{Casson-Gordon Knot Slicing Obstructions} There are a number of
formulations of how Casson-Gordon invariants provide obstructions to
slicing.  For our purposes the following fairly simple statement will
suffice.  Let $p$ be an odd prime and let $H_p$ denote 
 the $p$-primary summand of $H_1(M_K)$.  Again, we let $\beta$ denote the
torsion linking form.

\beginth{Theorem} If $K$ is slice there is a subgroup (or {\it
metabolizer}) 
$L_p \subset H_p$ with $|L_p|^2 = |H_p|$, $\beta(L_p,L_p) = 0$, and
$\tau(K,\chi) = 0$ for all
$\chi$ vanishing on $L_p$.\endth       

\section{The Main Theorem: p = 3} We now have the required material to
prove Theorem $0.1$.    To simplify notation, for any abelian group $A$,
let $A_p$ denote the $p$-primary summand of $A$.

\beginth{Theorem} If $|H_1(M_K)| = 3m$, where gcd($3,m$) = 1, then
$K$ is of infinite order in $\calc_1$.\endth

\noindent{\bf Proof.}  We have that $H_1(M_K)_3$ is isomorphic to
$\zz_3$,  generated by an element
$x$ with $\beta(x,x) = \pm({{1} \over {3}}) \in \qq/\zz$. Suppose now that
$K$ is of order $d$ in $\calc_1$.  Any nontrivial  metabolizing element 
for  
$H_1(M_{dK})_3
\cong (\zz_3)^d$ is of the form
$(x_i)_{i=1 \ldots d}$, where $x_i = \pm x$ for $r$ values of $i$ and is
$0$ otherwise.  Hence, Theorem $5.1$ yields that
$\tau(\#_{d}M_K,(\chi_{x_i})_{i=1 \ldots d})= 0$ with  exactly $r$ of the
$x_i = \pm1$ and  the rest $0$.  Here $\chi_{y}$ denotes the character
given by linking with $y$.

Now, applying Gilmer's additivity theorem and taking signatures we have
that
$r\sigma(K,\chi_1) = 0$.  (We have used that $\sigma(M_K,0) = 0$ and 
$\sigma(M_K,\chi_{-1}) = \sigma(M_K,\chi_1)$.)  It follows of course that
$\sigma(K,\chi_1) = 0$, contradicting Theorem 4.3.

\section{The Case p = 7}

The case of $p = 7$ introduces an issue we want to consider before dealing
with the general case.  One interesting observation about the following
argument is that the result does not follow from knowing the the vanishing
of the Casson-Gordon invariants for a spanning set of metabolizing
elements, or even their multiples. This demonstrates the very nonlinear
property of these invariants and  it is the first application we know of
in which that nonlinearity plays such an essential role.   

\beginth{Theorem} If $|H_1(M_K)| = 7m$, where gcd($7,m$) = 1, then
$K$ is not of order 4 in $\calc_1$.\endth

\vskip2ex

\noindent{\bf Proof }

 For such a knot K, a metabolizer for
$H_1(M_{4K})_7 \cong (\zz_7)^4$ can be seen to be generated by a pair of
elements $< (1, 0, 2,3), (0,1,-3,2)>$.  (There are other possibilities
differing only in order and sign from this one.)  Denoting by $\chi_a$ the
$\zz_7$-character that takes value $a$ on fixed generator of $H_1(M_K)_7$,
we find from either of these metabolizing vectors that $\tau(M_K,\chi_1) +
\tau(M_K,\chi_2) + \tau(M_K,\chi_3) = 0$. However, adding the generators
we see that the metabolizer must also contain the vector $(1,1,6,5)$ and
its multiples, $(2,2,5,3)$ and $(3,3,4,1)$.  From this we get the
relations $3\tau(M_K,\chi_1) +
\tau(M_K,\chi_2) = 0$, $3\tau(M_K,\chi_2) +
\tau(M_K,\chi_3) = 0$, and $3\tau(M_K,\chi_3) +
\tau(M_K,\chi_1) = 0$.  Combining these one finds that $28\tau(M_K,\chi_1)
= 0$ again contradicting Theorem 4.3.

\section{The General Case}

We now prove Theorem 0.1.

\beginth{Theorem}  If $|H_1(M_K)| = pm$ with $p$ a prime congruent to 3 
mod 4 and gcd(p, m) = 1, then $K$ is of infinite order in
$\calc_1$.\endth
\noindent{\bf Proof.}  Suppose that $dK$ is slice.  The existence of a
$Z_p$-metabolizer implies that $d$ is a multiple of 4. (The linking form of
$H_1(M_K)$ represents an element of order 4 in the Witt group of $Z_p$
linking forms.)  We begin by setting up some formalism to simplify the
sort of linear algebra that appeared in the previous section.  The example
below illustrates the notation we develop next.

Any metabolizing vector for the linking form on $H_1(M_K,\zz_p)$ (in the
subgroup $L_p$ given by Theorem 5.1) can be written as
$x = (x_i)_{i=1\ldots d} \in (\zz_p)^d$.  The condition that a
corresponding Casson-Gordon invariant vanishes yields $\sum_{x_i \ne 0}
\tau(M_K,\chi_{x_i}) = 0$.  Now the $x_i$ are in the cyclic group of
nonzero elements in $\zz_p$.  Denoting a generator for this group by $g$,
each nonzero $x_i$ corresponds to
$g^{\alpha_i}$ for some $\alpha_i$.  If we introduce further shorthand,
setting $t^{\alpha_i} = \tau(M_K,\chi_{x_i})$, we find that each
metabolizing vector leads to a relation $\sum_{x_i \ne 0} t^{\alpha_i} =
0$.  Note that at this point the symbol $t^{\alpha}$ does not represent a
power of any element ``$t$'', it is purely symbolic.  However it does
permit us to view
 the relations as being elements in the ring
$\zz[\zz_{p-1}]$.  Furthermore, since $\tau_{x_i} = \tau_{p-x_i}$, we have
that
$t^j = t^{j + (p-1) /2}$.  (Recall that $g^{(p-1)/{2}} = -1$.)  Hence, we
can view the relations as sitting in $\zz[\zz_{q}]$, where $q = (p-1)/2$.

Suppose that the metabolizing vector $x$ corresponds to the relation
$f = 0$, where $f$ is represented by an element in $\zz[\zz_q]$.  Then a
calculation shows that 
$ax$ corresponds to the relation $t^\alpha f$ where $g^\alpha = a$.  Hence
it follows that the relations between Casson-Gordon invariants generated
by a given element
$x
\in L_p$ and its multiples forms an ideal in $\zz[\zz_q]$ generated by the
polynomial $f$.  Before applying this to complete the proof, we should
pause for an example.

\vskip3ex

\noindent{\bf Example} Consider the metabolizing vector $x = (2,3,15,16)$
in
$(\zz_{19})^4$.  The nonzero elements of $\zz_{19}$ are generated by 2,
and we have $2 = 2^1$, $3 = 2^{13}$, $15 = 2^{11}$, and $16 = 2^4$.  Hence
in the notation just given, the vanishing of the corresponding
Casson-Gordon invariant can be written as $t^1 + t^{13} + t^{11} + t^4 =
0$.  Here we are in $\zz[\zz_{18}]$.  Switching to  $\zz[\zz_{9}]$ we have
that $t^1 + t^{4} + t^2 + t^4 = 0$.  (Notice that $\chi_{15}$ and $\chi_4$
yield the same Casson-Gordon invariant, and that $\chi_{15}$ corresponds
to $t^{11}$ while
$\chi_4$ corresponds to $t^2$ (since $2^2 = 4$) and $t^{11} = t^2 \in
\zz[\zz_{9}]$.)

Now consider the metabolizing vector $5x = (10, 15, 18, 4)$.  Since $10 =
2^{17}, 15 = 2^{11}, 18= 2^9$, and $4 = 2^2$, all mod 19, the
corresponding relation in $\zz[\zz_{18}]$ is $t^{17} + t^{11} + t^9 + t^2 =
0$.  Reducing to $\zz[\zz_{9}]$ gives $t^8 +t^2 + 1 + t^2$.

Hence by multiplying $x$ by 5 we have gone from the equation   $t^1 +
t^{4} + t^2 + t^4 = 0$ to  $t^8 +t^2 + 1 + t^2 =0$.  Notice that the
second polynomial is obtained from the first by multiplication by $t^7$. 
Finally $5 = 2^{16}$ mod 19, and
$t^{16} = t^7 \in \zz[\zz_{9}]$. 

\vskip3ex

To return to the proof, we must analyze the possible metabolizers $L_p$ for
$(\zz_p)^{4k}$, where $d = 4k$.  Such a metabolizer must be generated by
$2k$ elements.  Applying the Gauss-Jordan algorithm to a basis for $L_p$,
and perhaps reordering, we find a generating set $\{v_i\}_{i=1 \ldots
2k}$  where the first
$2k$ components of $v_i$ are 0, except the $i$-component which is 1.
Summing this basis produces the element $(1,1,\ldots,1,a_1,\ldots ,
a_{2k}) \in L_p$ where the first $2k$ entries are 1 and the $a_i$ are
unknown.

The corresponding relation is of the form $f =2k +
\sum_{i=1}^{k'}t^{\alpha_i} = 0$.  (The sum may not contain $2k$ terms if
any of the $a_i = 0$; hence $k'$ is less than or equal to $2k$.)  We next
show that the ideal generated by
$f$ in
$\zz[\zz_q]$ contains a nonzero integer.  This will follow from the fact
that
$f$ and $t^q -1$ are relatively prime, which will be the case unless $f$
vanishes at some
$q$-root of unity,
$\omega$; however, by considering norms and the triangle inequality we see
that this will be the case only if $k' = 2k$ and
$\omega^{a_i} = $ --1 for all $i$.  But since $q$ is odd, no power of
$\omega$ can equal --1.  

Since we now have that $f$ and $t^q -1$ are relatively prime, it follows
that with
$\qq$ coefficients (so that we are working over a PID) there is a
polynomial $g$ satisfying
$gf = 1$ mod ($t^q -1$).  Clearing denominators we find that for some
integral polynomial
$h$, $hf = n$ mod ($t^q -1$) for some positive integer $n$.

The proof of the theorem is  concluded by observing that we now have the
relation corresponding to $n \in \zz[\zz_q]$.  That is, $n\tau(M_K,\chi_1)
= 0$.  As before, this would imply that $\sigma(M_K,\chi_1) = 0$,
contradicting Theorem 4.3

\section{Corollaries}
\subsection{\bf Low crossing number knots}  Based on the work of Levine,
Morita [M] developed an algorithm to determine the order of a knot in the
algebraic concordance group using only its Alexander polynomial.  Based on
this, he enumerated all prime  knots of 10 or fewer crossings that are of
algebraic order 4.  There are eleven such knots, including $7_7$,
$9_{34}$, and nine 10 crossing knots.  Of these, seven have $H_1(M^3)$
satisfying our criteria for
$p=3$.  Three more satisfy the condition for $p=7$, and the last,
$10_{86}$, has   $H_1(M^3) = \zz_{83}$.  Hence, Corollary 0.2 follows.

\subsection{\bf Polynomial conditions}  A special case of Levine's results
states that a knot $K$ with $\Delta_K(t)$ quadratic is of finite order  if
$\Delta_K(t) = at^2 - (1+2a)t +a$ for some $a > 0$, and in that case it is
of order 4 if for some prime $p = 3$ mod 4, $\Delta_K(-1) = p^\alpha m$
with
$\alpha$ odd and gcd($p,m$) = 1. Our theorem applies only in the case that
$\alpha$ can be assumed to be 1.  However that is sufficient to give an
infinite family of examples, beginning with $\Delta_K(t) = 5t^2 - 11t +5$,
where $\Delta_K(-1) = (3)(7)$. 

\subsection{\bf Infinitely many linearly independent examples}  Knots
formed as twisted doubles of the unknot we among the first knots used to
construct algebraically slice knots that are not slice [CG1, CG2].  Jiang
[J] used these knots to demonstrate that the set of algebraically slice
knots contains a infinite set of knots that is  linearly independent in
concordance.  Here we demonstrate that twisted doubles also provide such
independent families  of knots that are of algebraic order 4.

To achieve independence within a family of examples, our theorem must be
extended somewhat.  Here is the statement we need. 

\beginth{Theorem}  If $|H_1(M_K)| = pm$ with $p$ a prime congruent to 3 
mod 4 and gcd$(p, m$) = 1, and if $J$ is any knot with $|H_1(M_J)| = q$,
where gcd($q,p$) = 1, then
$dK \# J $ is not slice for all nonzero integers $d$.\endth 

\noindent{\bf Proof.}  If we consider $\zz_p$ characters on the 2-fold
branched cover, the characters all vanish on $H_1(M_J)$ so by the
additivity of Casson-Gordon invariants we are reduced to considering the
character restricted to $H_1(M_{dK})$, which places us in the setting of
the proof of Theorem 0.1 in Section 8.

\vskip3ex To apply this, let $K_n$ denote the (--$n$)-twisted double of
the unknot, with
$n > 0$.  Then
$\Delta_{K_n}(t) = nt^2 - (1 + 2n)t + n$, and $H_1(M_{K_n}) = \zz_{4n
+1}$.  To pick an appropriate set of these knots, let
$\{ p_i \} $ be an enumeration of the primes that are congruent to 3 mod
4.  Let
$n_i = (p_{2i-1}p_{2i} - 1)/4$. Then the previous theorem quickly yields
the following.

\beginth{Corollary} The subset of the set of twisted doubles of the unknot
given by
$\{ K_{n_i} \} $, is a linearly independent set in the concordance group
and consists only of knots of algebraic order 4.\endth

\vfill
\eject

\vskip8ex
\centerline{\magbf References}
\vskip4ex

\item{[CG1]} A. Casson and C. Gordon {\sl Cobordism of classical knots,} in
``A la recherche de la Topologie perdue'', ed. by Guillou and Marin, 
Progress in Mathematics, Volume 62, 1986.

\item{[CG2]} A. Casson and C. Gordon {\sl On slice knots in dimension
three}, in Proc. Symp. Pure Math. 32 (1978), 39-54.

\item{[CF]} P.E. Conner and E. E.  Floyd, {\sl Differentiable periodic
maps}, Ergebnisse der Mathematik und ihrer Grenzgebiete, N. F., Band 33
Academic Press Inc., Publishers, New York; Springer-Verlag,
Berlin-Gšttingen-Heidelberg 1964.

\item{[CM]} D. Coray and F. Michel {\sl Knot cobordism and
amphicheirality}, Comm. Math. Helv., 58 (1983), 601-616.

\item{[FM]} R. Fox and J. Milnor, {\sl Singularities of $2$-spheres in
$4$-space and cobordism of knots}, Osaka J. Math. 3 (1966) 257--267. 

\item{[FQ]} M. Freedman and F. Quinn {\sl Topology of 4-manifolds},
Princeton Mathematical Series, 39. Princeton University Press, Princeton,
NJ, 1990.  

\item{[Gi]} P. Gilmer,  {\sl Slice knots in $S\sp{3}$},  Quart. J. Math.
Oxford Ser. (2) 34 (1983), no. 135, 305--322.

\item{[G]} C. Gordon, {\sl Problems}, in Knot Theory, ed. J.-C. Hausmann, 
Springer Lecture Notes no. 685, 1977.

\item{[J]} B. Jiang {\sl A simple proof that the concordance group of
algebraically slice knots is infinitely generated}, Proc. Amer. Math. Soc.
83 (1981), 189-192.

\item{[K1]} R. Kirby {\sl Problems in low dimensional manifold theory}, in
Algebraic and Geometric Topology (Stanford, 1976), vol 32, part II of Proc.
Sympos. Pure Math., 273--312.

\item{[K2]} R. Kirby {\sl Problems in low dimensional manifold theory}, in
Geometric Topology, AMS/IP Studies in Advanced Mathematics, ed. W. Kazez,
1997.

\item {[L1]} J. Levine {\sl Knot cobordism groups in codimension two}, 
Comment. Math. Helv. 44 (1969), 229--244.

 \item {[L2]} J. Levine {\sl Invariants of knot cobordism}, Invent. Math. 8
(1969), 98--110.

\item{[M]} T. Morita, {\sl Orders of knots in the algebraic knot cobordism
group}, Osaka J. Math. 25 (1988), no. 4, 859--864.   
   
\item{[T]} M. Toshiyuki {\sl On null equivalent knots}, Osaka J. Math. 11
(1959), no. 4, 95--113.

\vskip10ex

Department of Mathematics 

Indiana University

Bloomington, IN   47405

livingst@ucs.indiana.edu

\vskip3ex Department of Mathematics

University of Nevada

Reno, NV 89557

naik@unr.edu

\end